\documentclass[a4paper,11pt]{article} %Article class, provides some preliminary layout
\usepackage[utf8x]{inputenc} %gives us the symbols which are allowed to be used in the tex-code
\usepackage[english]{babel} %correct way to break words in english

\usepackage{amsmath} %Math enviroments, like align
\usepackage{hyperref} %References, citations, urls etc. can now be clicked 

\usepackage{amsfonts} %Font provided by AMS

\usepackage[parfill]{parskip} %To have no automated indendations after multiple linebreaks
\def    \date  {{ \number\day.\number\month.\number\year }}

\newcommand*{\doi}[1]{\url{https://doi.org/#1}} %DOI can be clicked and links to source

\newcommand {\rel} {{\mathbb R}}

\newcommand {\nat} {{\mathbb N}}

\newtheorem{proposition}{Proposition}[section]
\newtheorem{theorem}{Theorem}[section]
\newtheorem{lemma}[proposition]{Lemma}

\renewcommand{\theequation}{\mbox{\arabic{section}.\arabic{equation}}}

\newcommand{\nc}{\newcommand}
\newcommand{\renc}{\renewcommand}

\newcommand{\bcite}[1] {\cite{#1}}

\newcommand {\pr} {\bf}

\newcommand{\proof} {   \begin{flushright}
                        $\Box$%///
                        \end{flushright}
                }

\newcommand{\partproof} {       \begin{flushright}
                                //
                                \end{flushright}
                }

\newcommand{\defin} { \hspace*{\fill} $\Box$ }

\newcounter{claim}

\def    \Ln     {{ {\cal L}^n }}

\def	\d	{{ \ {\rm d} }}

\nc{\energ}[1]  {{ e_{#1} }}

\def	\Cn	{{ C }}
\def    \Cn     {{ C_{n} }}

\nc{\Ci}[1]	{{ C_{{#1}} }}
\renc{\Ci}[1]	{{ C_{n,{#1}} }}

\def	\poin	{{ poin }}

\nc{\gpo}[1]	{{ g_{\poin}^{#1} }}
\nc{\gpu}[1]	{{ g_{\poin,{#1}} }}

\nc{\gpmo}[1]	{{ g_{\poin,m}^{#1} }}
\nc{\tgpmo}[1]	{{ \tilde g_{\poin,m}^{#1} }}

\nc{\gplo}[1]	{{ g_{\poin,\lambda}^{#1} }}

\nc{\conf}[1]   {{ [{#1}] }}
\def	\lapm	{{ \ell^-_\varepsilon }}
\def	\lapp	{{ \ell^+_\varepsilon }}
\def	\lappm	{{ \ell^\pm_\varepsilon }}

\def	\Del	{{ M_\varepsilon }}

\newcommand {\winkel} %
{\mbox{$<$ \hskip -9pt \raisebox{1.1pt}{$\scriptscriptstyle )$}} \hskip 0.5ex}

\nc{\doub}[1]{{ \ddot{#1} }}
\nc{\dd}{ \begin{displaymath} }
\nc{\df}{ \end{displaymath} }
\nc{\dcd}{ \begin{displaymath} \begin{array}{c}}
\nc{\dcf}{ \end{array} \end{displaymath} }
\nc{\ee}{ \begin{equation} }
\nc{\ef}{ \end{equation} }
\nc{\ad}{ \begin{array}{c} }
\nc{\af}{ \end{array} }

\nc{\add}[1]    {}

\begin{document}

\begin{flushright}
\date
%submitted to
\end{flushright}

\begin{center}
{\LARGE \bf Positivity for the clamped plate equation under high tension}
\\ \ \\
Sascha Eichmann,
Reiner M. Sch\"atzle \\
Fachbereich Mathematik der
Eberhard-Karls-Universit\"at T\"ubingen, \\
Auf der Morgenstelle 10,
%Geb\"aude C, 5 A 40,
D-72076 T\"ubingen, Germany, \\
email: \href{mailto:sascha.eichmann@math.uni-tuebingen.de}{sascha.eichmann@math.uni-tuebingen.de}, \\
\href{mailto:schaetz@everest.mathematik.uni-tuebingen.de}{schaetz@everest.mathematik.uni-tuebingen.de} \\

\end{center}
\vspace{1cm}

\begin{quote}

{\bf Abstract:}
In this article we consider positivity issues for the clamped plate equation with high tension $\gamma>0$. This equation is given by $\Delta^2u - \gamma\Delta u=f$ under clamped boundary conditions. Here we show, that given a positive $f$, i.e. upwards pushing, we find a $\gamma_0>0$ such that for all $\gamma\geq \gamma_0$ the bending $u$ is indeed positive. This $\gamma_0$ only depends on the domain and the ratio of the $L^1$ and $L^\infty$ norm of $f$.\\ 
In contrast to a recent result by Cassani and Tarsia, our approach is valid in all dimensions.
\end{quote}

\vspace{1cm}

%\tableofcontents

%%%%%

\setcounter{equation}{0}

\section{Introduction} \label{intro}

The Boggio-Hadamard conjecture states, that for a given convex, open, bounded set $\Omega\subseteq\rel^n$, an $f\geq 0$ and outer normal $\nu$ of $\Omega$, a solution $u$ to
\begin{equation}
\label{eq:1_1}
 \left\{\begin{array}{cc}\Delta^2 u=f& \mbox{ in }\Omega,\\ u=\partial_\nu u=0& \mbox{ on }\partial\Omega.\end{array}\right.
\end{equation}
is nonnegative, i.e. $u\geq 0$ (cf. \cite{Hadamard1}, \cite{Hadamard2}). Problem (\ref{eq:1_1}) models the bending $u$ of a clamped plate $\Omega$ under a force $f$. Hence the problem can be restated as:
\begin{displaymath}
 \mbox{Does upward pushing yield upward bending?}
\end{displaymath}
The conjecture was substantiated by Boggio's explicit formula \cite{BoggioGreenFunction} (see also \cite[Lemma 2.27]{PolyHarmBoundValue} or \cite{HimmelBachelorThesis}) for the Greens function of problem (\ref{eq:1_1}) on the unit disc, because this function is positive, and furthermore by Almansi's result to calculate the Greens function for certain domains by the Greens function of the unit ball (see \cite{Almansi1899-1}). \\
Several other domains than the disc have been found on which such a positivity preserving property holds (see the references below). Hadamard himself claimed in \cite{Hadamard2} that such a property for all lima\c{c}ons is true, which turned out to be wrong in general, though some of them still possess this property (see \cite[Fig. 1.2]{PolyHarmBoundValue}). Remarkable is that such lima\c{c}ons are not convex. 
In \cite{GrunauRobert2010} Grunau and Robert showed, that positivity preserving is preserved under small regular perturbations of the domain in dimensions $n\geq 3$.
The conformal invariance of the problem was also successfully used to construct domains with a positivity preserving property by e.g. Dall'Acqua and Sweers in \cite{DallAcquaSweers2004} (see also the references therein for more informations on such domains).\\
On the other hand several counterexamples have been found by now. The first one was by Duffin on an infinite strip \cite{DuffinStrip} and shortly after Garabedian \cite{GarabedianEllipse} showed, that on an elongated ellipse the Greens function changes sign. By now even for uniform forces, i.e. $f\equiv1$, counterexamples have been found by Grunau and Sweers in \cite{GrunauSweersSignChangeUniform1} and \cite{GrunauSweersSignChangeUniform2}. We refer to \cite[§1.2]{PolyHarmBoundValue} for a comprehensive historical overview to this problem.\\
Instead of examining (\ref{eq:1_1}) for positivity, Cassani and Tarsia in \cite{CassaniTarsia} examined positivity issues for
\begin{equation}
\label{eq:1_2}
 \left\{\begin{array}{cc}\Delta^2 u - \gamma \Delta u=f& \mbox{ in }\Omega,\\ u=\partial_\nu u=0 &\mbox{ on }\partial\Omega.\end{array}\right.
\end{equation}
with $\gamma>0$ big enough. The basic motivation is that for $\gamma$ big enough, the influence of $\Delta u$ ($-\Delta u=f$ satisfies positivity preserving via the maximum principle) becomes stronger than that of $\Delta^2 u$. In more technical detail Cassani and Tarsia conjectured the existence of a $\gamma_0=\gamma_0(f,\Omega)\geq 0$, such that $u\geq0$ for all $\gamma\geq \gamma_0$ and provided a proof for dimensions $n=2,3$, smooth, bounded $\Omega$ and positive $f\in L^2(\Omega)$. 
In this article we provide a different approach, which is valid for all dimensions, see Theorem \ref{intro.theo}.\\
In dimension $n=1$ this positivity preserving property is true for all $\gamma>0$ independent of $f$. This was shown by Grunau in \cite{GrunauOneDimensionalCase2002} Proposition 1.\\
The parameter $\gamma$ is usually called tension, if it is positive. Several results concerning (\ref{eq:1_2}) have been achieved, which are usually concerned with vibrations of the plate, i.e. eigenvalue problems. Bickley gave some explicit caculations for the spectrum in the unit disc in \cite{Bickley} already in 1933. Hence the existing literature for these eigenvalue problems is quite vast and is still developing, see e.g. \cite{BuosoKennedy2021}, \cite{AntunesBuosoFreitag2019} and the references therein.\\
Other modifications for (\ref{eq:1_1}) concerning positivity issues are e.g. changing the boundary conditions to so called Steklov-conditions. This has been examined by e.g. Bucur and Gazzola in \cite{BucurGazzola2011}.\\
Different elliptic differential operators of higher order, their respective fundamental solutions and their sign close to a singularity have also been examined by Grunau, Romani and  Sweers in \cite{GrunauRomaniSweers2020} in a more systematic approach to understand better the loss of positivity preserving.

Instead of \eqref{eq:1_2} we examine the following boundary value problem for positivity preserving. This is obviously equivalent, but \eqref{intro.equ} yields the advantage, that the singularity of the equation is more prominent and hence yields easier access to necessary estimates.
\begin{equation} \label{intro.equ}
 \left\{\begin{array}{cc} \varepsilon^2\Delta^2 u_\varepsilon - \Delta u_\varepsilon = f&\mbox{ in }\Omega,\\
         u_\varepsilon=\partial_\nu u_\varepsilon=0&\mbox{ on } \partial\Omega.
        \end{array}\right.
\end{equation}
Here $\Omega\subset\subset \rel^n$, $\partial\Omega\in C^4$, $f\in L^\infty(\Omega)$ and $\varepsilon>0$,
and the solution $u_\varepsilon \in W^{4,p}(\Omega) \mbox{ for all } 1 < p < \infty$ (see e.g. \cite[Corollary 2.21]{PolyHarmBoundValue} and the references therein for existence, regularity and uniqueness to \eqref{intro.equ}).

\begin{theorem} \label{intro.theo}

For connected $\Omega\subset\subset \rel^n$, $\partial\Omega\in C^4$, $f\in L^\infty(\Omega)$, $\tau > 0, f \geq 0$ with
\begin{equation} \label{intro.theo.ass}
	\int f \d \Ln
	\geq \tau \parallel f \parallel_{L^\infty(\Omega)} > 0,
\end{equation}
there exists $\varepsilon_0 = \varepsilon_0(\Omega,\tau) > 0$ such that
\begin{equation} \label{intro.theo.conc}
	u_\varepsilon > 0
	\quad \mbox{in } \Omega
	\mbox{ for all } 0 < \varepsilon \leq \varepsilon_0.
\end{equation}
%\defin
\end{theorem}
Please note, that we do not have any restrictions on the dimension, i.e. $n\in\nat$ arbitrary. Furthermore our method yields, that $\varepsilon_0$ does not depend on $f$ directly, only on $\tau$.

The strategy of the proof is as follows:
The limiting problem of (\ref{intro.equ}) is
\begin{displaymath}
\left\{\begin{array}{cc}
		- \Delta u = f \geq 0, \not\equiv 0 & \mbox{ in } \Omega, \\
		u = 0 & \mbox{ on } \partial \Omega,
\end{array}\right.
\end{displaymath}
which admits a maximum principle,
and establishes positivity of $u_\varepsilon \mbox{ on any } \Omega' \subset \subset \Omega
\mbox{ for } \varepsilon$ small.

We proceed by contradiction and assume Theorem \ref{intro.theo} is false. Hence for every $\varepsilon>0$ we find a nonnegative $f_\varepsilon\in L^\infty(\Omega)$ satisfying \eqref{intro.theo.ass}, such that $u_\varepsilon$ is not positive in $\Omega$.

Then we examine a blow up of our solutions $u_\varepsilon$, which is weighted by the supremum of the modulus of the laplacian at the boundary, i.e. $\sup_{\partial \Omega} \varepsilon^2|\Delta u_\varepsilon|$. 
After a careful analysis (see sections \S \ref{prelim} and \S \ref{lap}), we can show that this blow up converges in a suitable sense to a solution of $\Delta^2 u-\Delta u=0$ on the half-space with Dirichlet boundary conditions (see section \S \ref{blow}). 
With a uniqueness result shown in appendix \ref{half}, we explicitly calculate this limit
and obtain positivity of the laplacian of $u_\varepsilon$ on the boundary for $\varepsilon$ small.
This is crucial, as in the presence of Dirichlet boundary conditions in (\ref{intro.equ})
the laplacian is the second normal derivative of $u_\varepsilon$ on the boundary,
and therefore positivity of the laplacian on the boundary gives positivity of $u_\varepsilon$
close to the boundary, see section \S \ref{blow}. \\
%This enables us to obtain informations on the blow up sequence itself, which in turn will yield our desired positivity preserving property at the end of section \S \ref{blow}.\\
Similar strategies of examining a blow up to the half space and using explicit formulas have been employed by Grunau and Robert in \cite{GrunauRobert2010} and Grunau, Robert and Sweers in \cite{GrunauRobertSweers2011} to show lower bounds for the Greens function of a polyharmonic operator. This method was later refined by Pulst in his PhD-thesis \cite{Pulst} to also obtain such estimates, if non-constant lower order terms are present. If variable coefficients in the principal part of the operator are given by a power of a second order elliptic linear operator, such estimates were found by the same method by Grunau in \cite{GrunauGreenEstimate2021}.\\
Our blow-up strategy needs careful estimates for the singular problem \eqref{intro.equ}. Estimates for these kinds of problems have a long history, see e.g. \cite{LFrank}, \cite{WGreenlee}, \cite{MGuegnon}, and \cite{Najman}. We are not aware of any specific estimates, which would help in our specific situation. For this reason and for the sake of completeness we derive them here.

\setcounter{equation}{0}

\section{Preliminary estimates} \label{prelim}
We proceed by contradiction and assume we find nonnegative $f_\varepsilon\in L^\infty(\Omega)$ for $\varepsilon\downarrow 0$ satisfying \eqref{intro.theo.ass}, such that the solution $u_\varepsilon$ of \eqref{intro.equ} is sign-changing. 
As (\ref{intro.equ}) is homogeneous of degree one,
we may assume by scaling
\begin{equation} \label{prelim.ass}
	0 \leq f_\varepsilon \leq 1\quad\mbox{ and }\quad
	\int f_\varepsilon \d \Ln \geq \tau > 0. \\
\end{equation}
By the Banach-Alaoglu Theorem we get after passing to subsequence and relabeling $f_\varepsilon \rightarrow f
\mbox{ weakly}^* \mbox{ in } L^\infty(\Omega)$ with
\begin{equation} \label{prelim.lim}
	0 \leq f \leq 1\quad \mbox{ and }\quad
	\int f \d \Ln \geq \tau > 0.
\end{equation}
In particular we have $f \not\equiv 0$.

The limiting problem of (\ref{intro.equ}) is thought to be the second-order boundary-value problem
\begin{equation} \label{prelim.equ}
\begin{array}{c}
	- \Delta u = f
	\quad \mbox{in } \Omega, \\
	u = 0
	\quad \mbox{on } \partial \Omega.
\end{array}
\end{equation}
We also consider
\begin{equation} \label{prelim.equ-eps}
\begin{array}{c}
	- \Delta u_{0,\varepsilon} = f_\varepsilon
	\quad \mbox{in } \Omega, \\
	u_{0,\varepsilon} = 0
	\quad \mbox{on } \partial \Omega.
\end{array}
\end{equation}
These two problems admit some important estimates:
\begin{proposition}
 \label{prelim.boundEstLim}
 For $u$ and $u_{0,\varepsilon}$ in \eqref{prelim.equ} rsp. \eqref{prelim.equ-eps} we find a constant $c_0=c_0(\Omega,\tau)>0$, such that
 \begin{equation} \label{prelim.u-normal}
	-\partial_\nu u, -\partial_\nu u_{0,\varepsilon} \geq c_0 > 0
	\quad \mbox{on } \partial \Omega
\end{equation}
and
\begin{equation} \label{prelim.u-low}
	u, u_{0,\varepsilon} \geq c_0 d(.,\partial \Omega)
	\quad \mbox{on } \Omega
\end{equation}
for $\varepsilon>0$ small enough.
\end{proposition}
{\pr Proof:}
Both problems admit by standard elliptic theory, see \bcite{gil.tru} Theorem 9.15,
unique solutions $u \mbox{ respectively } u_{0,\varepsilon}
\in W^{2,p}(\Omega) \hookrightarrow C^{1,\alpha}(\Omega)
\mbox{ for all } 1 < p < \infty \mbox{ and } 1 - (n/p) > \alpha > 0$ with
\begin{equation} \label{prelim.u-bound}
	\parallel u, u_{0,\varepsilon} \parallel_{W^{2,p}(\Omega) \cap C^{1,\alpha}(\Omega)}
	\leq C(\Omega,p,\alpha) < \infty
	\quad \forall 1 < p < \infty, 0 < \alpha < 1.
\end{equation}
In particular the set of all $u, u_{0,\varepsilon}
\mbox{ for any } f, f_\varepsilon$ with (\ref{prelim.lim}) and (\ref{prelim.ass})
is compact in $C^1(\overline{\Omega})$, and in particular
\begin{equation} \label{prelim.u-conv}
	u_{0,\varepsilon} \rightarrow u
	\quad \mbox{strongly in } C^{1,\alpha}(\Omega),
	\forall 0 < \alpha < 1,
\end{equation}
and weakly in $W^{2,p}(\Omega) \mbox{ for all } 1 < p < \infty$.

As $f, f_\varepsilon \geq 0, f, f_\varepsilon \not\equiv 0 \mbox{ and } \Omega$ is connected,
we get by the strong maximum principle, see \bcite{gil.tru} Theorem 8.19, that
\begin{equation} \label{prelim.u-pos}
	u, u_{0,\varepsilon} > 0 \mbox{ in } \Omega.
\end{equation}
By Hopf's maximum principle,
see \bcite{gil.tru} Lemma 3.4, and by compactness of the $u, u_{0,\varepsilon} \mbox{ in } C^1(\overline{\Omega})$ we find a constant $c_0 = c_0(\Omega,\tau) > 0$ such that \eqref{prelim.u-normal} and \eqref{prelim.u-low} both hold for $\varepsilon>0$ small enough.
\proof
We put $v_\varepsilon := u_\varepsilon - u_{0,\varepsilon}
\in W^{2,p}(\Omega) \cap C^{1,\alpha}(\Omega)
\mbox{ with } v_\varepsilon = 0 \mbox{ on } \partial \Omega$.
Subtracting in (\ref{intro.equ}), we see
\begin{equation} \label{prelim.harm}
	\Delta (\varepsilon^2 \Delta u_\varepsilon - v_\varepsilon) = 0
	\quad \mbox{in } \Omega,
\end{equation}
that is $\varepsilon^2 \Delta u_\varepsilon - v_\varepsilon$
is harmonic in $\Omega$.

Equation (\ref{prelim.equ}) is indeed the limiting problem of (\ref{intro.equ})
in the sense of the following proposition.

\begin{proposition} \label{prelim.conv}

For $u_\varepsilon, u, v_\varepsilon = u_\varepsilon - u_{0,\varepsilon}$
as in (\ref{intro.equ}), (\ref{prelim.ass}), (\ref{prelim.equ-eps}),
we have
\begin{equation} \label{prelim.conv-ass}
\begin{array}{c}
	u_\varepsilon \rightarrow u, v_\varepsilon \rightarrow 0
	\mbox{ strongly in } W^{1,2}_0(\Omega), \\
	\varepsilon \Delta u_\varepsilon \rightarrow 0 \mbox{ strongly in } L^2(\Omega).
\end{array}
\end{equation}
\end{proposition}
{\pr Proof:} \\
Multiplying (\ref{prelim.harm}) by $v_\varepsilon$ and integrating by parts, we get
\begin{displaymath}
	0 = \int \limits_\Omega \Delta (\varepsilon^2 \Delta u_\varepsilon - v_\varepsilon)
	v_\varepsilon \d \Ln =
\end{displaymath}
\begin{displaymath}
	= \int \limits_\Omega \varepsilon^2 (\Delta^2 u_\varepsilon) u_\varepsilon \d \Ln
	- \int \limits_\Omega \varepsilon^2 (\Delta^2 u_\varepsilon) u_{0,\varepsilon} \d \Ln
	+ \int \limits_\Omega |\nabla v_\varepsilon|^2 \d \Ln.
\end{displaymath}
Replacing $\varepsilon^2 \Delta^2 u_\varepsilon \mbox{ by } \Delta v_\varepsilon$
in the second term with (\ref{prelim.harm}), we continue
\begin{displaymath}
	\int \limits_\Omega \varepsilon^2 |\Delta u_\varepsilon|^2 \d \Ln
	+ \int \limits_\Omega |\nabla v_\varepsilon|^2 \d \Ln
	= \int \limits_\Omega (\Delta v_\varepsilon) u_{0,\varepsilon} \d \Ln =
\end{displaymath}
\begin{equation} \label{prelim.aux}
	= -\int \limits_\Omega  \nabla v_\varepsilon \nabla u_{0,\varepsilon} \d \Ln
	\leq \parallel \nabla v_\varepsilon \parallel_{L^2(\Omega)}
	\ \parallel \nabla u_{0,\varepsilon} \parallel_{L^2(\Omega)},
\end{equation}
%\begin{displaymath}
%	= \int \limits_\Omega |\nabla u_{0,\varepsilon}|^2 \d \Ln
%	- \int \limits_\Omega \nabla u_\varepsilon \nabla u_{0,\varepsilon} \d \Ln.
%\end{displaymath}
in particular $\parallel \nabla v_\varepsilon \parallel_{L^2(\Omega)}
\leq \parallel \nabla u_{0,\varepsilon} \parallel_{L^2(\Omega)}
\leq C(\Omega)$ by (\ref{prelim.u-bound}), hence
\begin{displaymath}
	\int \limits_\Omega \varepsilon^2 |\Delta u_\varepsilon|^2 \d \Ln
	+ \int \limits_\Omega |\nabla v_\varepsilon|^2 \d \Ln
	\leq \int \limits_\Omega |\nabla u_{0,\varepsilon}|^2 \d \Ln
	\leq C(\Omega).
\end{displaymath}
Passing to a subsequence, we get
\begin{displaymath}
\begin{array}{c}
	v_\varepsilon \rightarrow v \mbox{ weakly in } W^{1,2}_0(\Omega), \\
	\varepsilon^2 \Delta u_\varepsilon \rightarrow 0 \mbox{ strongly in } L^2(\Omega).
\end{array}
\end{displaymath}
Multiplying (\ref{prelim.harm}) by some $\eta \in C^\infty_0(\Omega)$
and passing to a subsequence, we get
\begin{displaymath}
	0 = \int \limits_\Omega (\varepsilon^2 \Delta^2 u_\varepsilon - \Delta v_\varepsilon) \eta \d \Ln
	= \int \limits_\Omega \varepsilon^2 \Delta u_\varepsilon \cdot \Delta \eta \d \Ln
	+ \int \limits_\Omega \nabla v_\varepsilon \nabla \eta \d \Ln
	\rightarrow \int \limits_\Omega \nabla v \nabla \eta \d \Ln,
\end{displaymath}
and $v \in W^{1,2}_0(\Omega)$ is harmonic in $\Omega$,
hence $v = 0 \mbox{ and } u_\varepsilon \rightarrow u
\mbox{ weakly in } W^{1,2}_0(\Omega)$.
Returning to (\ref{prelim.aux}), we improve now to
\begin{displaymath}
	\int \limits_\Omega \varepsilon^2 |\Delta u_\varepsilon|^2 \d \Ln
	+ \int \limits_\Omega |\nabla v_\varepsilon|^2 \d \Ln
	= -\int \limits_\Omega \nabla v_\varepsilon \nabla u_{0,\varepsilon} \d \Ln
	\rightarrow 0,
\end{displaymath}
which is (\ref{prelim.conv-ass}).
\proof
The following proposition shows that the laplacian cannot be bounded throughout $\Omega$.

\begin{proposition} \label{prelim.int}

For $u_\varepsilon, \tau$
as in (\ref{intro.equ}), (\ref{prelim.ass}), (\ref{prelim.lim}),
we have for any $\Omega' \subset \subset \Omega$ that
\begin{displaymath}
	\liminf \limits_{\varepsilon \rightarrow 0}
	\int \limits_{\Omega \setminus \Omega'} (\Delta u_\varepsilon)_+ \d \Ln
	\geq \tau > 0.
\end{displaymath}
\end{proposition}
{\pr Proof:} \\
We see for any $\eta \in C^\infty_0(\Omega)$
with the previous Proposition \ref{prelim.conv}
\begin{displaymath}
	\int (\Delta u_\varepsilon) \eta \d \Ln
	= \int u_\varepsilon \Delta \eta \d \Ln
	\rightarrow \int u \Delta \eta \d \Ln
	= \int (\Delta u) \eta \d \Ln
	= -\int f \eta \d \Ln
\end{displaymath}
and the homogeneous boundary conditions in (\ref{intro.equ}) that
\begin{displaymath}
	\int \Delta u_\varepsilon (1 - \eta) \d \Ln
	= \int \Delta u_\varepsilon \d \Ln
	- \int (\Delta u_\varepsilon) \eta \d \Ln =
\end{displaymath}
\begin{displaymath}
	= \int \limits_{\partial \Omega} \partial_\nu u_\varepsilon \d area_{\partial \Omega}
	- \int (\Delta u_\varepsilon) \eta \d \Ln
	\rightarrow \int f \eta \d \Ln.
\end{displaymath}
Choosing $\eta \in C^\infty_0(\Omega) \mbox{ with } 0 \leq \eta \leq 1
\mbox{ and } \eta \equiv 1 \mbox{ in } \Omega'$, we get
\begin{displaymath}
	\int \limits_{\Omega \setminus \Omega'} (\Delta u_\varepsilon)_+ \d \Ln
	\geq \int (\Delta u_\varepsilon) (1 - \eta) \d \Ln
	\rightarrow \int f \eta \d \Ln.
\end{displaymath}
Letting $\eta \nearrow \chi_\Omega$, we get from (\ref{prelim.lim}) that
\begin{displaymath}
	\int f \eta \d \Ln \rightarrow \int f \d \Ln \geq \tau > 0,
\end{displaymath}
and the proposition follows.
\proof

%%%%%
%%%%%

\setcounter{equation}{0}

\section{The laplacian on the boundary} \label{lap}

In this section, we investigate the values of the laplacian on the boundary and put
\begin{equation} \label{lap.lap-boundary-inter}
	\lapm := \min \limits_{\partial \Omega} \varepsilon^2 \Delta u_\varepsilon
	\leq \max \limits_{\partial \Omega} \varepsilon^2 \Delta u_\varepsilon =: \lapp.
\end{equation}
With subscripts $\pm$ we denote the positive respectively negative part, i.e.
\begin{equation}
 \label{lap.lap-pos-neg-part}
 \lappm_{,+}:=\max\left(0,\lappm\right),\ \lappm_{,-}:=\max\left(0, - \lappm\right). 
\end{equation}
Furthermore we set
\begin{equation} \label{lap.lap-boundary-sup}
	\Del
	:= \parallel \varepsilon^2 \Delta u_\varepsilon \parallel_{L^\infty(\partial \Omega)}
	= \max( \lapp_{,+} , \lapm_{,-} ).
\end{equation}
The quantity $M_\varepsilon$ will be crucial throughout the exposition. Our goal is to show, that it has the same asymptotic as $\varepsilon$ itself, i.e. we find constants $c_0,C>0$ such that
\begin{equation*}
 c_0 \varepsilon\leq M_\varepsilon\leq C\varepsilon.
\end{equation*}
for $\varepsilon>0$ small.
A first step in this direction is Proposition \ref{lap.asymp}, which will later be improved to our desired result in Proposition \ref{blow.asymp} and \eqref{blow.fini.asymp}.

With the maximum principle, we get the following estimates.

\begin{proposition} \label{lap.max}

For $u_\varepsilon, f_\varepsilon, u_{0,\varepsilon},
v_\varepsilon = u_\varepsilon - u_{0,\varepsilon}, \lappm$
as in (\ref{intro.equ}), (\ref{prelim.ass}), (\ref{prelim.equ-eps}),
(\ref{lap.lap-boundary-inter}),
we have
\begin{equation} \label{lap.max.harm}
	\lapm \leq \varepsilon^2 \Delta u_\varepsilon - v_\varepsilon \leq \lapp
	\quad \mbox{in } \Omega,
\end{equation}
\begin{equation} \label{lap.max.v}
	-\lapp_{,+} - \varepsilon^2 \parallel f_\varepsilon \parallel_{L^\infty(\Omega)}
	\leq v_\varepsilon
	\leq \lapm_{,-}
	\quad \mbox{in } \Omega,
\end{equation}
\begin{equation} \label{lap.max.u-low}
	u_\varepsilon > -\lapp_{,+}
	\quad \mbox{in } \Omega.
\end{equation}
\end{proposition}
{\pr Proof:} \\
As $v_\varepsilon = u_\varepsilon - u_{0,\varepsilon} = 0 \mbox{ on } \partial \Omega$
by (\ref{intro.equ}) and (\ref{prelim.equ-eps}),
we get (\ref{lap.max.harm}) from (\ref{prelim.harm}).

Adding (\ref{prelim.equ-eps}), we see
\begin{displaymath}
	\varepsilon^2 \Delta v_\varepsilon - v_\varepsilon
	= -\varepsilon^2 \Delta u_{0,\varepsilon} + \varepsilon^2 \Delta u_\varepsilon - v_\varepsilon
	\leq \varepsilon^2 f_\varepsilon + \lapp
	\quad \mbox{in } \Omega,
\end{displaymath}
in particular
\begin{displaymath}
	\varepsilon^2 \Delta v_\varepsilon < 0
	\quad \mbox{in } [v_\varepsilon < -\varepsilon^2 f_\varepsilon - \lapp_{,+}] =: \Omega_0.
\end{displaymath}
Since $v_\varepsilon \geq -\varepsilon^2 f_\varepsilon - \lapp_{,+}
\mbox{ on } \partial \Omega_0$,
as $v_\varepsilon = 0 \mbox{ on } \partial \Omega$ by above,
we get from the mean-value estimate for superharmonic functions
or by Alexandroff's maximum principle, as $v_\varepsilon \in W^{2,n}(\Omega)$,
see \bcite{gil.tru} Theorem 9.1,
that $v_\varepsilon \geq -\varepsilon^2 f_\varepsilon - \lapp_{,+}
\mbox{ in } \Omega_0$,
hence $\Omega_0 = \emptyset$,
and the left estimate in (\ref{lap.max.v}) follows.
The right estimate is obtained by symmetry observing that $f_\varepsilon \geq 0$.

Next for $x \in \overline{\Omega} \mbox{ with }
u_\varepsilon(x) = \min_{\overline{\Omega}} u_\varepsilon$ and assuming that this minimum is negative,
we see $x \in \Omega$, as $u_\varepsilon = 0 \mbox{ on } \partial \Omega$ by (\ref{intro.equ}),
hence $\Delta u_\varepsilon(x) \geq 0$,
as $u_\varepsilon \in W^{4,p}(\Omega) \hookrightarrow C^2(\Omega)
\mbox{ for } 2 - (n/p) > 0$.
Then we get with (\ref{lap.max.harm}) and (\ref{prelim.u-pos}) that
\begin{displaymath}
	-u_\varepsilon(x)
	\leq \varepsilon^2 \Delta u_\varepsilon(x) - v_\varepsilon(x) - u_{0,\varepsilon}(x)
	< \lapp,
\end{displaymath}
which is (\ref{lap.max.u-low}).
\proof
Using the fourth order equation, we get estimates for the laplacian.

\begin{proposition} \label{lap.lap}

For $u_\varepsilon, f_\varepsilon, u_{0,\varepsilon},
v_\varepsilon = u_\varepsilon - u_{0,\varepsilon}, \lappm$
as in (\ref{intro.equ}), (\ref{prelim.ass}), (\ref{prelim.equ-eps}),
(\ref{lap.lap-boundary-inter}),
we have
\begin{equation} \label{lap.lap.u}
	-\lapm_{,-} - \varepsilon^2 \parallel f_\varepsilon \parallel_{L^\infty(\Omega)}
	\leq \varepsilon^2 \Delta u_\varepsilon
	\leq \lapp_{,+}
	\quad \mbox{in } \Omega.
\end{equation}
\end{proposition}
{\pr Proof:} \\
We have with (\ref{intro.equ}) that
\begin{displaymath}
	\varepsilon^2 \Delta^2 u_\varepsilon - \Delta u_\varepsilon = f_\varepsilon
	\quad \mbox{in } \Omega
\end{displaymath}
and get
\begin{displaymath}
	\varepsilon^2 \Delta (\Delta u_\varepsilon) < 0
	\quad \mbox{in } [\Delta u_\varepsilon
	< -\parallel f_{\varepsilon,+} \parallel_{L^\infty(\Omega)}
	- \varepsilon^{-2} \lapm_{,-}] =: \Omega_0.
\end{displaymath}
Since $\Delta u_\varepsilon \geq -\parallel f_{\varepsilon,+} \parallel_{L^\infty(\Omega)}
- \varepsilon^{-2} \lapm_{,-} \mbox{ on } \partial \Omega_0$,
as $\Delta u_\varepsilon \geq \varepsilon^{-2} \lapm
\geq -\varepsilon^{-2} \lapm_{,-} \mbox{ on } \partial \Omega$ with (\ref{lap.lap-boundary-inter}),
we get from the mean-value estimate for superharmonic functions
or by Alexandroff's maximum principle, as $\Delta u_\varepsilon \in W^{2,n}(\Omega)$,
see \bcite{gil.tru} Theorem 9.1,
that $\Delta u_\varepsilon \geq -\parallel f_{\varepsilon,+} \parallel_{L^\infty(\Omega)}
- \varepsilon^{-2} \lapm_{,-} \mbox{ in } \Omega_0$,
hence $\Omega_0 = \emptyset$,
and the left estimate in (\ref{lap.lap.u}) follows.
The right estimate is obtained by symmetry observing that $f_{\varepsilon} \geq 0$.
\proof
Here we can give a preliminary asymptotic estimate for $\Del$.
Actually we will improve this asymptotic later
in Proposition \ref{blow.asymp} and (\ref{blow.fini.asymp}).
Anyway we present this estimate at this stage
to get more compact bounds already now.

\begin{proposition} \label{lap.asymp}

For $u_\varepsilon, \Del, \lappm$
as in (\ref{intro.equ}), (\ref{prelim.ass}),
(\ref{lap.lap-boundary-inter}), (\ref{lap.lap-boundary-sup}),
we have
\begin{displaymath}
	\varepsilon^{-2} \Del \geq \varepsilon^{-2} \lapp_{,+} \rightarrow \infty.
\end{displaymath}
in particular $\Del \geq \lapp = \lapp_{,+} > 0 \mbox{ for } \varepsilon$ small
depending on $\Omega \mbox{ and } \tau$.
\end{proposition}
{\pr Proof:} \\
Combining Proposition \ref{lap.lap} (\ref{lap.lap.u}) and Proposition \ref{prelim.int},
we get for any $\Omega' \subset \subset \Omega$ that
\begin{displaymath}
	\Ln(\Omega \setminus \Omega') \liminf \limits_{\varepsilon \rightarrow 0}
	\varepsilon^{-2} \lapp_{,+}
	\geq \liminf \limits_{\varepsilon \rightarrow 0}
	\int \limits_{\Omega \setminus \Omega'} (\Delta u_\varepsilon)_+ \d \Ln
	\geq \tau > 0,
\end{displaymath}
hence, as $\Ln(\Omega \setminus \Omega')$ can be made arbitrarily small, that
\begin{displaymath}
	\varepsilon^{-2} \lapp_{,+} \rightarrow \infty,
\end{displaymath}
which yields the assertion.
\proof
With the above asymptotic,
we can already bound $v_\varepsilon$,
and we can prove that $u$ is positive on large parts of $\Omega$.

\begin{proposition} \label{lap.v}

For $u_\varepsilon, u_{0,\varepsilon},
v_\varepsilon = u_\varepsilon - u_{0,\varepsilon}, \Del$
as in (\ref{intro.equ}), (\ref{prelim.ass}), (\ref{prelim.lim}), (\ref{prelim.equ-eps}),
(\ref{lap.lap-boundary-sup}),
we have
\begin{displaymath}
	\limsup \limits_{\varepsilon \rightarrow 0}
	\parallel \Del^{-1} v_\varepsilon \parallel_{L^\infty(\Omega)} \leq 1.
\end{displaymath}
\end{proposition}
{\pr Proof:} \\
Combining Proposition \ref{lap.max} (\ref{lap.max.v}) and Proposition \ref{lap.asymp},
we get observing (\ref{prelim.ass}) that
\begin{displaymath}
	\limsup \limits_{\varepsilon \rightarrow 0}
	\parallel \Del^{-1} v_\varepsilon \parallel_{L^\infty(\Omega)}
	\leq \limsup \limits_{\varepsilon \rightarrow 0}
	\Del^{-1} (\Del + \varepsilon^2 \parallel f_\varepsilon \parallel_{L^\infty(\Omega)})
	\leq 1.
\end{displaymath}
\proof

\begin{proposition} \label{lap.u-pos-far}

For $u_\varepsilon, \lappm$
as in (\ref{intro.equ}), (\ref{prelim.ass}), (\ref{lap.lap-boundary-inter}),
we have
\begin{equation} \label{lap.u-pos-far.ass}
	u_\varepsilon > 0
	\quad \mbox{in } [d(.,\partial \Omega) \geq C \lapp] \cap \Omega
\end{equation}
for some $C = C(\Omega,\tau) < \infty \mbox{ and } \varepsilon$ small
depending on $\Omega \mbox{ and } \tau$.
\end{proposition}
{\pr Proof:} \\
From (\ref{prelim.u-low}), we see for $x \in \Omega
\mbox{ with } c_0 d(x,\partial \Omega)
> \lapp_{,+} + \varepsilon^2 \parallel f_\varepsilon \parallel_{L^\infty(\Omega)}$
by (\ref{lap.max.v}) that
\begin{displaymath}
	u_\varepsilon(x)
	= v_\varepsilon(x) + u(x)
	\geq -\lapp_{,+} - \varepsilon^2 \parallel f_\varepsilon \parallel_{L^\infty(\Omega)}
	+ c_0 d(x,\partial \Omega) > 0.
\end{displaymath}
By Proposition \ref{lap.asymp} and (\ref{prelim.ass}) clearly
\begin{displaymath}
	2 \lapp \geq \lapp_{,+} + \varepsilon^2 \parallel f_\varepsilon \parallel_{L^\infty(\Omega)}
\end{displaymath}
for $\varepsilon$ small,
and (\ref{lap.u-pos-far.ass}) follows for $C = 2 c_0^{-1} < \infty$.
\proof

%%%%%
%%%%%

\setcounter{equation}{0}

\section{Blow up} \label{blow}

In this section, we consider a blow up of our solutions $u_\varepsilon$
by translating and rescaling with $x_{0,\varepsilon} \in \rel^n$. We will have to choose $x_{0,\varepsilon}$ differently in different steps in the proof of Theorem \ref{intro.theo}. %For this please note, that the results of Propositions \ref{blow.conv}, \ref{blow.asymp}, \ref{blow.conv-0} and \ref{blow.lim} do not depend on a fixated $x_{0,\varepsilon}$, but it is allowed to be varied. 
In e.g. the proof of Claim \ref{blow.fini.asymp-c} in the proof of Theorem \ref{intro.theo} below, we choose $x_{0,\varepsilon}\in \partial \Omega$
 such that $\max_{x\in\partial \Omega}|\varepsilon^2 \Delta u_\varepsilon(x)|$ is attained,
 while in the proof of Claim \ref{blow.fini.lap-pos-c} below, we choose $x_{0,\varepsilon}\in\partial \Omega$ to attain $\lapm$.
This will not cause a problem, because the constants yielded by these claims do not depend on $\varepsilon$.
We put
\begin{equation} \label{blow.def}
	\tilde u_\varepsilon(x) := u_\varepsilon(x_{0,\varepsilon} + \varepsilon x)
	\quad \mbox{for } x \in \tilde \Omega_\varepsilon := \varepsilon^{-1} (\Omega - x_{0,\varepsilon}).
\end{equation}
Then
\begin{equation} \label{blow.lap}
	\nabla \tilde u_\varepsilon = \varepsilon
	\nabla u_\varepsilon(x_{0,\varepsilon} + \varepsilon .),
	\quad \Delta \tilde u_\varepsilon
	= \varepsilon^2 \Delta u_\varepsilon(x_{0,\varepsilon} + \varepsilon .)
\end{equation}
and by (\ref{intro.equ}) that
\begin{equation} \label{blow.equ}
\begin{array}{c}
	\Delta^2 \tilde u_\varepsilon - \Delta \tilde u_\varepsilon
	= \varepsilon^2 f_\varepsilon(x_{0,\varepsilon} + \varepsilon .) =: \tilde f_\varepsilon
	\quad \mbox{on } \tilde \Omega_\varepsilon, \\
	\tilde u_\varepsilon, \partial_{\nu_{\tilde \Omega_\varepsilon}} \tilde u_\varepsilon = 0
	\quad \mbox{on } \partial \tilde \Omega_\varepsilon
\end{array}
\end{equation}
and with (\ref{prelim.ass}) and Proposition \ref{lap.asymp} that
\begin{equation} \label{blow.equ-f}
	\parallel \tilde f_\varepsilon \parallel_{L^\infty(\tilde \Omega_\varepsilon)}
	= \parallel f_\varepsilon \parallel_{L^\infty(\Omega)} \varepsilon^2 = o(1) \Del
	\quad \mbox{for } \varepsilon \rightarrow 0,
\end{equation}
where $o(1) \rightarrow 0$ depending on $\Omega \mbox{ and } \tau$.
We extend $u_\varepsilon \mbox{ respectively } \tilde u_\varepsilon$
by putting $0 \mbox{ outside } \Omega \mbox{ respectively outside } \tilde \Omega_\varepsilon$.
By the homogeneous boundary conditions in (\ref{intro.equ}) and (\ref{blow.equ}),
we see $u_\varepsilon, \tilde u_\varepsilon \in W^{2,2}_{loc}(\rel^n)$.

We also have to stretch $\tilde u_\varepsilon$ to get a nontrivial limit,
and it turns out
that reaching bounded values of the laplacian of $\tilde u_\varepsilon$
on the boundary is the right measure for stretching.

\begin{proposition} \label{blow.conv}

For $\tilde u_\varepsilon, \tilde f_\varepsilon, u_\varepsilon, f_\varepsilon, \Del$
as in (\ref{blow.def}), (\ref{intro.equ}), (\ref{prelim.ass}), (\ref{lap.lap-boundary-sup}),
and $x_{0,\varepsilon} \in \partial \Omega$,
we get for any subsequence with after rotating $\tilde \Omega_\varepsilon$ such that
\begin{equation} \label{blow.conv.normal}
	\nu_\Omega(x_{0,\varepsilon}) \rightarrow -e_n
\end{equation}
after passing to a subsequence
\begin{equation} \label{blow.conv.conv}
\begin{array}{c}
	\Del^{-1} \tilde u_\varepsilon \rightarrow \tilde u_\infty
	\quad \mbox{weakly in } W^{2,2}(B_R(0)) \mbox{ for all } R > 0, \\
	\Del^{-1} \tilde u_\varepsilon \rightarrow \tilde u_\infty
	\quad \mbox{weakly in } W^{4,p}(B_R(0) \cap \tilde \Omega_\varepsilon)
	\mbox{ for all } R > 0, 1 < p < \infty
\end{array}
\end{equation}
as $\varepsilon \rightarrow 0$
after flattening the boundary of $\partial \tilde \Omega_\varepsilon$.
Further
\begin{equation} \label{blow.conv.u-lim}
\begin{array}{c}
	\Delta^2 \tilde u_\infty - \Delta \tilde u_\infty = 0
	\quad \mbox{in } \rel^n_+, \\
	\tilde u_\infty, \partial_n \tilde u_\infty = 0
	\quad \mbox{on } \rel^{n-1} \times \{ 0 \}, \\
	\tilde u_\infty = 0
	\quad \mbox{in } \rel^n \setminus \rel^n_+, \\
	\tilde u_\infty \geq -1,
	\quad |\Delta \tilde u_\infty| \leq 1
	\quad \mbox{in } \rel^n. \\
\end{array}
\end{equation}
\end{proposition}
{\pr Proof:} \\
We get from Proposition \ref{lap.max} (\ref{lap.max.u-low}) that
\begin{equation} \label{blow.conv.u-low}
	\tilde u_\varepsilon \geq -\Del
	\quad \mbox{in } \tilde \Omega_\varepsilon
\end{equation}
and also outside $\tilde \Omega_\varepsilon$,
as $\tilde u_\varepsilon = 0$ there.
Next by Proposition \ref{lap.lap} (\ref{lap.lap.u}),
Proposition \ref{lap.asymp}, (\ref{prelim.ass}) and (\ref{blow.lap}) that
\begin{equation} \label{blow.conv.lap-infty}
	\parallel \Delta \tilde u_\varepsilon \parallel_{L^\infty(\tilde \Omega_\varepsilon)}
	= \parallel \varepsilon^2 \Delta u_\varepsilon
	\parallel_{L^\infty(\Omega)}
	\leq \Del + \varepsilon^2 \parallel f_\varepsilon \parallel_{L^\infty(\Omega)}
	\leq (1 + o(1)) \Del.
\end{equation}
Then $\tilde u_\varepsilon + \Del \geq 0 \mbox{ in } \rel^n$,
and we can apply the Harnack-inequality, see \bcite{gil.tru} Theorems 8.17 and 8.18,
and get observing that $\tilde u_\varepsilon(0) = 0$,
as $x_{0,\varepsilon} \in \partial \Omega$, that
\begin{displaymath}
	\sup \limits_{B_R(0)} (\tilde u_\varepsilon + \Del)
	\leq \Cn (\tilde u_\varepsilon + \Del)(0)
	+ \Ci{R} \parallel \Delta (\tilde u_\varepsilon + \Del) \parallel_{L^\infty(B_{2R}(0))}
	\leq \Ci{R} \Del,
\end{displaymath}
hence
\begin{equation} \label{blow.conv.u-infty}
	\parallel \tilde u_\varepsilon \parallel_{L^\infty(B_R(0))}
	\leq \Ci{R} \Del
\end{equation}
for $\varepsilon$ small depending on $\Omega \mbox{ and } \tau$.

Next by Friedrichs's Theorem in the interior, see \bcite{gil.tru} Theorem 8.8, \bcite{gil.tru} Exercise 8.2,
(\ref{blow.conv.lap-infty}) and (\ref{blow.conv.u-infty}) that
\begin{equation} \label{blow.conv.w22}
	\parallel \tilde u_\varepsilon \parallel_{W^{2,2}(B_R(0))}
	\leq \Ci{R} \Big( \parallel \Delta \tilde u_\varepsilon \parallel_{L^2(B_{2R}(0))}
	+ \parallel \tilde u_\varepsilon \parallel_{L^2(B_{2R}(0))} \Big)
	\leq \Ci{R} \Del
\end{equation}
for $\varepsilon$ small depending on $\Omega \mbox{ and } \tau$,
which yields the first convergence in (\ref{blow.conv.conv})
after passing to a subsequence.

Proceeding from (\ref{blow.equ}),
we get from fourth order $L^p-$estimates,
see \bcite{adn1}, \bcite{adn2} \S 10,
after flattening the boundary of $\partial \tilde \Omega_\varepsilon$,
as $\partial \Omega \in C^4$,
with (\ref{blow.equ-f}) and (\ref{blow.conv.u-infty}) that
\begin{displaymath}
	\parallel \tilde u_\varepsilon \parallel_{W^{4,p}(B_R(0) \cap \tilde \Omega_\varepsilon)} \leq
\end{displaymath}
\begin{displaymath}
	\leq \Ci{R,p} \Big( \parallel \tilde f_\varepsilon
	\parallel_{L^\infty(B_{2R}(0) \cap \tilde \Omega_\varepsilon)}
	+ \parallel \tilde u_\varepsilon
	\parallel_{L^\infty(B_{2R}(0) \cap \tilde \Omega_\varepsilon)} \Big) \leq
\end{displaymath}
\begin{displaymath}
	\leq \Ci{R,p} \Del
	\quad \forall R > 0, 1 < p < \infty
\end{displaymath}
and $\varepsilon$ small depending on $\Omega \mbox{ and } \tau$.
After passing to a subsequence,
we obtain with (\ref{blow.conv.normal}) the second convergence in (\ref{blow.conv.conv}).

Finally (\ref{blow.conv.u-lim}) follows from (\ref{blow.equ}), (\ref{blow.equ-f}),
(\ref{blow.conv.u-low}), (\ref{blow.conv.lap-infty}), 
when recalling that $\tilde u_\varepsilon = 0
\mbox{ in } \rel^n \setminus \Omega_\varepsilon$.
\proof
Actually by fourth order higher order $L^p-$estimates,
see \bcite{adn1}, \bcite{adn2} \S 10,
we get that the blow up $\tilde u_\infty$ is smooth on $\overline{\rel^n_+}$.

Now we are able to give a lower bound for $\Del$
which improves the asymptotic in Proposition \ref{lap.asymp}.

\begin{proposition} \label{blow.asymp}

For $u_\varepsilon, \Del$
as in (\ref{intro.equ}), (\ref{prelim.ass}),
(\ref{lap.lap-boundary-sup}),
we have
\begin{displaymath}
	\Del \geq c_0 \varepsilon
\end{displaymath}
for some $c_0 = c_0(\Omega,\tau) > 0 \mbox{ and } \varepsilon$ small
depending on $\Omega \mbox{ and } \tau$.
\end{proposition}
{\pr Proof:} \\
We see for $x_\varepsilon \in \Omega
\mbox{ with } |x_\varepsilon - x_{0,\varepsilon}|
= d(x_\varepsilon,\partial \Omega) = \varepsilon
\mbox{ for } \varepsilon$ small
by (\ref{prelim.u-low}), Proposition \ref{lap.v}
and by the local boundedness of $\tilde u_\varepsilon$
in Proposition \ref{blow.conv}, or more directly by (\ref{blow.conv.u-infty}),
for $\tilde x_\varepsilon := (x_\varepsilon - x_{0,\varepsilon}) / \varepsilon
\in \overline{B_1(0)}$ that
\begin{displaymath}
	c_0 \varepsilon
	\leq u_{0,\varepsilon}(x_\varepsilon)
	= u_\varepsilon(x_\varepsilon) - v_\varepsilon(x_\varepsilon)
	= \tilde u_\varepsilon(\tilde x_\varepsilon) - v_\varepsilon(x_\varepsilon) \leq
\end{displaymath}
\begin{displaymath}
	\leq \parallel \tilde u_\varepsilon \parallel_{L^\infty(B_1(0))}
	+ \parallel v_\varepsilon \parallel_{L^\infty(\Omega)}
	\leq (\Ci{1} + o(1)) \Del,
\end{displaymath}
hence
\begin{displaymath}
	\Del \geq c_0 \varepsilon
\end{displaymath}
for $\varepsilon$ small depending on $\Omega \mbox{ and } \tau$.
\proof
The blow up for $u_{0,\varepsilon}$ is rather elementary by the strong convergence
in $C^{1,\alpha}(\Omega)$ in (\ref{prelim.u-conv}).
As in (\ref{blow.def}), we put
\begin{equation} \label{blow.def-0}
	\tilde u_{0,\varepsilon}(x) := u_{0,\varepsilon}(x_{0,\varepsilon} + \varepsilon x)
	\quad \mbox{for } x \in \tilde \Omega_\varepsilon := \varepsilon^{-1} (\Omega - x_{0,\varepsilon}).
\end{equation}

\begin{proposition} \label{blow.conv-0}

For $\tilde u_{0,\varepsilon}, u_\varepsilon, u_{0,\varepsilon}, f_\varepsilon, \Del$
as in (\ref{blow.def-0}), (\ref{intro.equ}), (\ref{prelim.ass}), (\ref{lap.lap-boundary-sup}),
and with (\ref{blow.conv.normal}),
we have after passing to a subsequence such that
\begin{equation} \label{blow.conv-0.beta}
	\beta \leftarrow |\nabla u_{0,\varepsilon}(x_{0,\varepsilon})| \varepsilon/\Del
	\geq c_0 \varepsilon/\Del \geq 0
\end{equation}
for some $c_0 = c_0(\Omega,\tau) > 0$
for the linear function $\tilde u_{0,\infty}: (y,t) \rightarrow \beta t$ that
\begin{displaymath}
	\Del^{-1} \tilde u_{0,\varepsilon} \rightarrow \tilde u_{0,\infty}
	\quad \mbox{uniformly in } B_R(0) \cap \tilde \Omega_\varepsilon
	\mbox{ for all } R > 0
\end{displaymath}
after flattening the boundary of $\partial \tilde \Omega_\varepsilon$.
Further
\begin{equation} \label{blow.conv-0.linear}
	|\tilde u_\infty - \tilde u_{0,\infty}| \leq 1
	\quad \mbox{in } \rel^n_+,
\end{equation}
in particular $\tilde u_\infty - \tilde u_{0,\infty} \in L^\infty(\rel^n_+)$.
\end{proposition}
{\pr Proof:} \\
From (\ref{prelim.u-bound}) and,
as $u_{0,\varepsilon} = 0 \mbox{ on } \partial \Omega$ by (\ref{prelim.equ-eps}),
we get by Taylor's expansion for any $x \in \rel^n_+
\mbox{ and any } 0 < \alpha < 1$ that
\begin{equation}\label{blow.conv-0.taylor}
	\Del^{-1} \tilde u_{0,\varepsilon}(x)
	= \Del^{-1} u_{0,\varepsilon}(x_{0,\varepsilon} + \varepsilon x)
	= \Del^{-1} \varepsilon \nabla u_{0,\varepsilon}(x_{0,\varepsilon}) \cdot x
	+ \Del^{-1} O_\alpha(\varepsilon^{1 + \alpha} |x|^{1 + \alpha}).
\end{equation}
As $u_{0,\varepsilon} = 0 \mbox{ on } \partial \Omega$ by (\ref{prelim.equ-eps}),
we get with (\ref{prelim.u-conv}), (\ref{blow.conv.normal}) and $\partial \Omega \in C^4$
after passing to a subsequence with $x_{0,\varepsilon}
\rightarrow x_0 \in \partial \Omega$ that
\begin{displaymath}
	\nabla u_{0,\varepsilon}(x_{0,\varepsilon})
	= \partial_{\nu_\Omega} u_{0,\varepsilon}(x_{0,\varepsilon}) \nu_\Omega(x_{0,\varepsilon})
	\rightarrow -\partial_{\nu_\Omega} u(x_0) e_n,
\end{displaymath}
in particular with (\ref{prelim.u-normal}) that
\begin{displaymath}
	|\nabla u_{0,\varepsilon}(x_{0,\varepsilon})|
	\rightarrow -\partial_{\nu_\Omega} u(x_0) > 0
\end{displaymath}
and $\nabla u_{0,\varepsilon}(x_{0,\varepsilon}) / |\nabla u_{0,\varepsilon}(x_{0,\varepsilon})|
\rightarrow e_n$. 
By Proposition \ref{blow.asymp} we extract a subsequence, such that $\varepsilon /M_\varepsilon$ converges for $\varepsilon\downarrow 0$. Then (\ref{prelim.u-normal}) yields
\begin{displaymath}
	\beta := \leftarrow
	|\nabla u_{0,\varepsilon}(x_{0,\varepsilon})| \varepsilon/\Del
	\geq c_0 \varepsilon /\Del \geq 0,
\end{displaymath}
which is (\ref{blow.conv-0.beta}). Furthermore Proposition \ref{blow.asymp} yields
\begin{displaymath}
	\varepsilon^{1 + \alpha} / \Del \rightarrow 0.
\end{displaymath}
Together, we get
\begin{displaymath}
	\Del^{-1} \varepsilon \nabla u_{0,\varepsilon}(x_{0,\varepsilon})
	\rightarrow \beta e_n,
\end{displaymath}
and the proposed convergence follows from the Taylor expansion \eqref{blow.conv-0.taylor}.

Further we get with Proposition \ref{lap.v} that
\begin{displaymath}
	|\tilde u_\infty - \tilde u_{0,\infty}|
	\leftarrow |\Del^{-1} \tilde u_\varepsilon - \Del^{-1} \tilde u_{0,\varepsilon}|
	\leq \limsup \limits_{\varepsilon \rightarrow 0}
	\parallel \Del^{-1} v_\varepsilon \parallel_{L^\infty(\Omega)} \leq 1,
\end{displaymath}
which is (\ref{blow.conv-0.linear}).
\proof
By our investigation of half space solutions in Appendix \ref{half},
Proposition \ref{half.u} applied to $\tilde u_\infty \mbox{ and } \tilde u_{0,\infty}$
with Proposition \ref{blow.conv} (\ref{blow.conv.u-lim})
and Proposition \ref{blow.conv-0} (\ref{blow.conv-0.linear})
determines $\tilde u_\infty$ uniquely as the one dimensional solution
and immediately yields the following Proposition.

\begin{proposition} \label{blow.lim}

For $\tilde u_\infty$ as in Proposition \ref{blow.conv}
and $\beta$ as in Proposition \ref{blow.conv-0},
we have
\begin{equation} \label{blow.lim.null}
	\tilde u_\infty \equiv 0,
	\quad \mbox{if } \beta = 0,
\end{equation}
\begin{equation} \label{blow.lim.lap}
	\Delta \tilde u_\infty \equiv \beta \geq 0
	\quad \mbox{in } \rel^{n-1} \times \{ 0 \}
\end{equation}
and
\begin{equation} \label{blow.lim.u-pos}
	\tilde u_\infty > 0
	\quad \mbox{in } \rel^n_+,
	\quad \mbox{if } \beta > 0.
\end{equation}
\defin
\end{proposition}
Now we are able to conclude the proof of Theorem \ref{intro.theo}.
\\ \ \\
{\pr Proof of Theorem \ref{intro.theo}:} \\
We consider $u_\varepsilon, f_\varepsilon, u_{0,\varepsilon}, \lappm, \Del$
as in (\ref{intro.equ}), (\ref{prelim.ass}), (\ref{prelim.equ-eps}),
(\ref{lap.lap-boundary-inter}), (\ref{lap.lap-boundary-sup})
and $\tilde u_\varepsilon, \tilde f_\varepsilon, \tilde u_{0,\varepsilon}$
as in (\ref{blow.def}), (\ref{blow.def-0})
with their blow ups $\tilde u_\infty, \tilde u_{0,\infty} \mbox{ and } \beta$
obtained in the Propositions \ref{blow.conv} and \ref{blow.conv-0}.
We prove various claims.

\setcounter{claim}{0}
\ \\
\refstepcounter{claim} \label{blow.fini.asymp-c}
{\pr Claim \arabic{claim}: \\ }
\begin{displaymath}
	\liminf \limits_{\varepsilon \rightarrow 0}
	\varepsilon / \Del > 0,
\end{displaymath}
hence with Proposition \ref{blow.asymp} that
\begin{equation} \label{blow.fini.asymp}
	c_0 \leq \Del / \varepsilon \leq C
\end{equation}
for some $c_0 = c_0(\Omega,\tau) > 0, C = C(\Omega,\tau) < \infty
\mbox{ and } \varepsilon$ small
depending on $\Omega \mbox{ and } \tau$. \\
{\pr Proof:} \\
If on contrary $\varepsilon / \Del \rightarrow 0$
for a subsequence $\varepsilon \rightarrow 0$,
then we get from Proposition \ref{blow.conv-0}
and (\ref{prelim.u-bound}) that $\beta = 0$,
hence with Proposition \ref{blow.lim} (\ref{blow.lim.null}) after passing to this subsequence we have
that $\tilde u_\infty \equiv 0$.

On the other hand choosing $x_{0,\varepsilon} \in \partial \Omega$ in such a way that
\begin{displaymath}
	|\varepsilon^2 \Delta u_\varepsilon(x_{0,\varepsilon})|
	= \parallel \varepsilon^2 \Delta u_\varepsilon \parallel_{L^\infty(\partial \Omega)}
	= \Del,
\end{displaymath}
we get from the convergence in Proposition \ref{blow.conv} (\ref{blow.conv.conv}) that
\begin{equation} \label{blow.fini.asymp-aux}
	1 = \Del^{-1} |\varepsilon^2 \Delta u_\varepsilon(x_{0,\varepsilon})|
	= |\Del^{-1} \Delta \tilde u_\varepsilon(0)|
	\rightarrow |\Delta \tilde u_\infty(0)|,
\end{equation}
hence $\Delta \tilde u_\infty(0) \neq 0 \mbox{ and } \tilde u_\infty \not\equiv 0$.
This is a contradiction, and the claim follows.
\partproof
\refstepcounter{claim} \label{blow.fini.lap-pos-c}
{\pr Claim \arabic{claim}: \\ }
\begin{equation} \label{blow.fini.lap-lim-pos}
	\liminf \limits_{\varepsilon \rightarrow 0} \Del^{-1} \lapm > 0,
\end{equation}
\begin{equation} \label{blow.fini.lap-pos}
	\Del^{-1} \Delta \tilde u_\varepsilon \geq c_1 > 0
	\quad \mbox{on } \partial \tilde \Omega_\varepsilon
\end{equation}
for some $c_1 = c_1(\Omega,\tau) > 0 \mbox{ and } \varepsilon$ small
depending on $\Omega \mbox{ and } \tau$
and
\begin{equation} \label{blow.fini.u-pos}
	\tilde u_\infty > 0
	\quad \mbox{in } \rel^n_+.
\end{equation}
{\pr Proof:} \\
(\ref{blow.fini.asymp}) implies with Proposition \ref{blow.conv-0} (\ref{blow.conv-0.beta}) that
\begin{displaymath}
	\beta
%	\geq \tilde \beta \liminf \limits_{\varepsilon \rightarrow 0} \varepsilon/\Del
	\geq c_0 \liminf \limits_{\varepsilon \rightarrow 0} \varepsilon/\Del
	> 0,
\end{displaymath}
which immediately gives (\ref{blow.fini.u-pos}) by Proposition \ref{blow.lim} (\ref{blow.lim.u-pos}).

Next we choose $x_{0,\varepsilon} \in \partial \Omega$ in such a way that
\begin{displaymath}
	\varepsilon^2 \Delta u_\varepsilon(x_{0,\varepsilon})
	= \min \limits_{\partial \Omega} \varepsilon^2 \Delta u_\varepsilon
	= \lapm
\end{displaymath}
and get as in (\ref{blow.fini.asymp-aux})
from the convergence in Proposition \ref{blow.conv} (\ref{blow.conv.conv})
and Proposition \ref{blow.lim} (\ref{blow.lim.lap}) that
\begin{displaymath}
	\Del^{-1} \lapm
	= \Del^{-1} \varepsilon^2 \Delta u_\varepsilon(x_{0,\varepsilon})
	= \Del^{-1} \Delta \tilde u_\varepsilon(0)
	\rightarrow \Delta \tilde u_\infty(0) = \beta > 0,
\end{displaymath}
and (\ref{blow.fini.lap-lim-pos}) follows.
Clearly (\ref{blow.fini.lap-lim-pos}) implies with (\ref{lap.lap-boundary-inter}) that
\begin{displaymath}
	\Del^{-1} \inf \limits_{\partial \tilde \Omega_\varepsilon} \Delta \tilde u_\varepsilon
	= \Del^{-1} \inf \limits_{\partial \Omega} \varepsilon^2 \Delta u_\varepsilon
	= \Del^{-1} \lapm \geq c_1 > 0
\end{displaymath}
for some $c_1 > 0 \mbox{ and } \varepsilon$ small, which is (\ref{blow.fini.lap-pos}).
\partproof
Further (\ref{blow.fini.lap-lim-pos}) implies
that $\lapp \geq \lapm > 0 \mbox{ and } \lapm_{,-} = 0$ for $\varepsilon$ small,
hence with (\ref{lap.lap-boundary-sup}) that
\begin{equation} \label{blow.fini.lap-max}
	\Del = \max( \lapp_{,+} , \lapm_{,-} ) = \lapp
\end{equation}
for $\varepsilon$ small depending on $\Omega \mbox{ and } \tau$.
\refstepcounter{claim} \label{blow.fini.pos-far-c}
\\ \ \\
{\pr Claim \arabic{claim}: \\ }
\begin{equation} \label{blow.fini.pos-far}
	u_\varepsilon > 0
	\quad \mbox{in } [d(.,\partial \Omega) \geq C \varepsilon] \cap \Omega
\end{equation}
for some $C = C(\Omega,\tau) < \infty \mbox{ and } \varepsilon$ small
depending on $\Omega \mbox{ and } \tau$. \\
{\pr Proof:} \\
This follows directly from Proposition \ref{lap.u-pos-far},
(\ref{blow.fini.asymp}) and (\ref{blow.fini.lap-max}).
\partproof
\refstepcounter{claim} \label{blow.fini.pos-near-c}
{\pr Claim \arabic{claim}: \\ }
\begin{equation} \label{blow.fini.pos-near}
	u_\varepsilon > 0
	\quad \mbox{in } [0 < d(.,\partial \Omega) \leq c_2 \varepsilon] \cap \Omega
\end{equation}
for some $c_2 = c_2(\Omega,\tau) > 0 \mbox{ and } \varepsilon$ small
depending on $\Omega \mbox{ and } \tau$. \\
{\pr Proof:} \\
By the homogeneous boundary conditions in (\ref{intro.equ}),
we have $\Delta u_\varepsilon(x_{0,\varepsilon}) = \partial_{\nu \nu} u_\varepsilon(x_{0,\varepsilon})$,
hence with (\ref{blow.fini.lap-pos}), Proposition \ref{blow.conv} (\ref{blow.conv.conv})
and the embedding $W^{4,p} \hookrightarrow C^3 \mbox{ for } 1 - (n/p) > 0$ that
\begin{displaymath}
	\tilde u_\varepsilon(-t \nu_{\tilde \Omega_\varepsilon}(0))
	\geq \frac{1}{2} t^2 \Delta \tilde u_\varepsilon(0) - \Cn \Del t^3
	\geq \frac{1}{2} t^2 \Del (c_1 - 2 \Cn t) > 0
\end{displaymath}
for $0 < t < c_1 / (2 \Cn) < 1 \mbox{ and } \varepsilon$ small.
As $x_{0,\varepsilon} \in \partial \Omega$ can be chosen arbitrarily,
we get for $c_2 = c_1 / (4 \Cn)$ that
\begin{displaymath}
	\tilde u_\varepsilon > 0
	\quad \mbox{in } [0 < d(.,\partial \tilde \Omega_\varepsilon) \leq c_2]
	\cap \tilde \Omega_\varepsilon,
\end{displaymath}
and the claim follows by rescaling in (\ref{blow.def}).
\partproof
\refstepcounter{claim} \label{blow.fini.pos-middle-c}
{\pr Claim \arabic{claim}: \\ }
\begin{equation} \label{blow.fini.pos-middle}
	u_\varepsilon > 0
	\quad \mbox{in } [c_2 \varepsilon \leq d(.,\partial \Omega) \leq C \varepsilon] \cap \Omega
\end{equation}
for $\varepsilon$ small
depending on $\Omega \mbox{ and } \tau$. \\
{\pr Proof:} \\
Here for any $x_\varepsilon \in \Omega
\mbox{ with } c_2 \varepsilon \leq d(x_\varepsilon,\partial \Omega) \leq C \varepsilon
\mbox{ for } \varepsilon$ small,
we select $x_{0,\varepsilon} \in \partial \Omega
\mbox{ with } d(x_\varepsilon,\partial \Omega) = |x_\varepsilon - x_{0,\varepsilon}|$
and get for $\tilde x_\varepsilon := \varepsilon^{-1} (x_\varepsilon - x_{0,\varepsilon})
\in \tilde \Omega_\varepsilon$ that
\begin{displaymath}
	c_2 \leq d(\tilde x_\varepsilon,\partial \tilde \Omega_\varepsilon)
	= |\tilde x_\varepsilon| \leq C.
\end{displaymath}
Passing to subsequence, we get $\tilde x_\varepsilon \rightarrow \tilde x
\mbox{ with } \tilde x \in \overline{\rel^n_+}$ and
\begin{displaymath}
	d(\tilde x , \rel^{n-1} \times \{0\}) \geq c_2 > 0,
\end{displaymath}
hence $\tilde x \in \rel^n_+$.
Then by Proposition \ref{blow.conv} (\ref{blow.conv.conv}) and (\ref{blow.fini.u-pos}) that
\begin{displaymath}
	\Del^{-1} u_\varepsilon(x_\varepsilon)
	= \Del^{-1} \tilde u_\varepsilon(\tilde x_\varepsilon)
	\rightarrow \tilde u_\infty(\tilde x) > 0,
\end{displaymath}
and we conclude that $u_\varepsilon(x_\varepsilon) > 0
\mbox{ for } \varepsilon$ small,
and the claim follows.
\partproof
Combining (\ref{blow.fini.pos-far}), (\ref{blow.fini.pos-near}) and (\ref{blow.fini.pos-middle}),
we get
\begin{displaymath}
	u_\varepsilon > 0
	\quad \mbox{in } \Omega
\end{displaymath}
for $\varepsilon$ small depending on $\Omega \mbox{ and } \tau$,
which proves Theorem \ref{intro.theo}.
\proof

%%%%%

%%%%%
\ \\
{\LARGE \bf Appendix}
%\\ \ \\
%In this appendix, we collect for the reader's convenience
%some results which are consequences or adaptions of standard
%results.

\begin{appendix}
\renewcommand{\theequation}{\mbox{\Alph{section}.\arabic{equation}}}

%%%%%

\setcounter{equation}{0}

\section{Uniqueness of a bi-laplace equation on the half space} \label{half}

In this section we show the following uniqueness theorem.
\begin{theorem}
\label{half.fourth}
 Let $v:\overline{\rel^{n}_+}\rightarrow \rel$ be a smooth solution to
 $$\left\{\begin{array}{c}
           \Delta^2 v - \Delta v=0 \mbox{ on }\rel^{n}_+,\\
           v=\partial_n v =0  \mbox{ on }\rel^{n-1}\times\{0\},\\
           v\in L^\infty(\rel^{n}_+).
          \end{array}\right.$$
Then $v= 0$.
%\defin
\end{theorem}
The proof is based on an energy type estimate and on the one-dimensional case. 
We start with the one dimensional case and show the following lemma:
\begin{lemma}
\label{A_2}
 Let $v\in C^4_{loc}([0,\infty[)$ satisfy
%\begin{displaymath}
%\begin{array}{cc}
   %\left\{\begin{array}{c}
   \begin{align*}
    \frac{d^4 v}{dt^4} - \frac{d^2v}{dt^2}=0& \mbox{ on }]0,\infty[,\\
           v(0)=\frac{d v}{dt}(0) =0,  \\
           v\in L^\infty(]0,\infty[).
   \end{align*}
          %\end{array}\right.
%\end{array}
%\end{displaymath}
Then $v=0$.
\end{lemma}
{\pr Proof:} \\
 Since the differential equation is ordinary and linear, the solution space of the equation itself is of dimension $4$. By inserting the following functions, we see that they constitute a basis of the solution space
 $$t\mapsto \cosh(t),\ \sinh(t),\ 1,\ t.$$
 The following functions are therefore a basis of the solution space including the initial conditions at $t=0$:
 $$v_1(t)=\cosh(t)-1,\ v_2(t)=\sinh(t)-t.$$
 Hence $v$ has to be of the form
 $$v(t)=A v_1(t) + Bv_2(t) = A(\cosh(t)-1) + B(\sinh(t)-t).$$
 Inserting the exponential function for $\cosh$ and $\sinh$ yields
\begin{align*}
v(t) =& A \left(\frac{e^t + e^{-t}}{2} -1 \right) + B \left(\frac{e^t - e^{-t}}{2} - t\right)\\ 
=& \frac{1}{2}(A+B) e^t + \frac{1}{2}(A-B)e^{-t} - A - Bt.
\end{align*}
Since the solution is bounded, we have $A=-B$. Again the boundedness then yields $B=0$, to rule out linear growth. Hence $v=0$.
\proof
For our next step we introduce a bit of notation for half spaces
$$B_R^+ := B_R(0) \cap \rel^{n}_+ \subseteq \rel^n,\ R>0.$$
Next we show an energy type estimate.
\begin{lemma}
 \label{A_3}
 Let $v$ satisfy
 $$\left\{\begin{array}{c}
           \Delta^2 v - \Delta v=0 \mbox{ on }\rel^{n}_+,\\
           v=\partial_n v =0  \mbox{ on }\rel^{n-1}\times\{0\}.
          \end{array}\right.$$
Then there exists a constant $C=C(n)>0$, such that for all $R>1$ we have
$$\int_{B_R^+} \left(|\Delta v|^2 + |\nabla v|^2\right)\, d\Ln \leq \frac{C}{R^2}\int_{B_{2R}^+}\left(|\nabla v|^2 + v^2\right)\, d\Ln.$$
\end{lemma}
{\pr Proof:} \\
 Let $\eta\in C^\infty_0(B_{2R}(0))$ with $0\leq \eta\leq 1$, $\eta=1$ on $B_R(0)$ and for any $k\geq1$
\begin{equation}
\label{eq:A_1}
 |D^k \eta|\leq \frac{C_k}{R^k}\chi_{B_{2R}(0)\setminus B_R(0)},
\end{equation}
i.e. $\eta$ is a cutoff function for the ball $B_R(0)$. Then $v \eta^4$ and its first derivative are zero on $\partial B_{2R}^+$. Therefore partial integration and the differential equation itself yield
\begin{align*}
0 =& \int_{B_{2R}^+} v\eta^4(\Delta v - \Delta^2v)\,d\Ln
= -\int_{B_{2R}^+}\Delta(v\eta^4)\Delta v \,d\Ln - \int_{B_{2R}^+}\nabla(v\eta^4)\nabla v\,d\Ln \\
=& - \int_{B_{2R}^+} \left(|\Delta v|^2\eta^4 + |\nabla v|^2\eta^4\right)\,d\Ln\\
&- \int_{B_{2R}^+}2 \left(\nabla v\nabla(\eta^4)\Delta v+ v\Delta(\eta^4)\Delta v + v\nabla v \nabla (\eta^4)\right)\,d\Ln\\ 
=& - \int_{B_{2R}^+} \left(|\Delta v|^2\eta^4 + |\nabla v|^2\eta^4\right)\,d\Ln\\
&-\int_{B_{2R}^+} \big(8\nabla v\nabla\eta \eta^3\Delta v+ 4v\Delta v\Delta \eta \eta^3 + 12v\Delta v|\nabla\eta|^2\eta^2+ 4 v\nabla v\nabla \eta \eta^3\big)\,d\Ln.
\end{align*}
By the previous identity and Young's inequality with an $\varepsilon>0$ we get
\begin{align*}
\int_{B_{2R}^+}& \left(|\Delta v|^2\eta^4 + |\nabla v|^2\eta^4\right) \,d\Ln\\
\leq \int_{B_{2R}^+}&\big(\varepsilon|\Delta v|^2\eta^4+ C_\varepsilon |\nabla v|^2|\nabla\eta|^2\eta^2\\
&+ \varepsilon|\Delta v|^2 \eta ^4+ C_\varepsilon v^2 |\Delta \eta|^2\eta^2\\
&+ \varepsilon|\Delta v|^2 \eta^4 + C_\varepsilon v^2|\nabla\eta|^4\\
&+ \varepsilon |\nabla v|^2 \eta ^4 + C_\varepsilon v^2|\nabla\eta|^2\eta^2\big) \,d\Ln.
\end{align*}
Now choosing $\varepsilon$ small enough, we can absorb the terms with $\varepsilon$ as a prefactor into the left hand side %and obtain with $0\leq \eta \leq 1$
\begin{displaymath}
	\frac{1}{2}\int_{B_{2R}^+} \left(|\Delta v|^2\eta^4 + |\nabla v|^2\eta^4\right)\,d\Ln
\end{displaymath}
\begin{displaymath}
	\leq C\int_{B_{2R}^+}\left(|\nabla v|^2|\nabla\eta|^2\eta^2
	+v^2 |\Delta \eta|^2\eta^2+v^2|\nabla\eta|^4+v^2|\nabla\eta|^2\eta^2\right)\,d\Ln.
\end{displaymath}
Together with $0\leq\eta\leq1$ the estimates on the derivatives of the cutoff function (\ref{eq:A_1}) yield
$$\int_{B_R^+}\left(|\Delta v|^2 + |\nabla v|^2\right)\,d\Ln\leq C\int_{B_{2R}^+}\left(|\nabla v|^2\frac{1}{R^2}+v^2 \frac{1}{R^4}+v^2\frac{1}{R^4}+v^2\frac{1}{R^2}\right)\,d\Ln.$$
Since we have choosen $R>1$, we obtain
$$\int_{B_R^+}\left(|\Delta v|^2 + |\nabla v|^2\right)\,d\Ln\leq \frac{C}{R^2}\int_{B_{2R}^+}\left(|\nabla v|^2+v^2\right)\,d\Ln.$$
\proof
Now we can show our main result Theorem \ref{half.fourth}, by iterating Lemma \ref{A_3}:
\\ \ \\
{\pr Proof of Theorem \ref{half.fourth}:} \\
 Since $v$ is bounded and the differential equation is linear and elliptic, we can employ Schauder-type estimates (see \cite[Thm. 6.2]{adn1} and \cite[Sec. 4]{DouglisNirenbergInteriorSchauder}) to obtain
 \begin{equation}
\label{eq:A_2}
 \|v\|_{C^{k}(\overline{\rel^{n}_+})}\leq C_k\|v\|_{L^\infty(\overline{\rel^{n}_+})}<\infty.
\end{equation}
Let $x=(y,t)\in\overline{\rel^{n}_+}$, such that $y\in\rel^{n-1}$ and $t\geq 0$. By $\partial_y^k v$ we denote any partial derivative  of $v$ of order $k$ only after horizontal directions, i.e. indices in $\{1,\ldots,n-1\}$. Then for every $k\in\nat$ the function $\partial^k_{y} v:\overline{\rel^n_+}\rightarrow\rel $ still solves the differential equation, satisfies the Dirichlet boundary conditions and by (\ref{eq:A_2}) is again bounded. Now we can iteratively apply Lemma \ref{A_3} for $k\geq\ell\in\nat$ and obtain
\begin{displaymath}
	\int_{B_R^+}|\partial^k_{y} v|^2\, d\Ln
	\leq \int_{B_R^+}|\nabla\partial^{k-1}_{y} v|^2\, d\Ln
\end{displaymath}
\begin{displaymath}
	\leq \frac{C}{R^2}\int_{B_{2R}^+}
	\left(|\nabla\partial^{k-1}_{y} v|^2 + |\partial^{k-1}_{y} v|^2\right)\, d\Ln
\end{displaymath}
\begin{displaymath}
	\leq \frac{C}{R^2}\int_{B_{2R}^+}\left(|\nabla\partial^{k-1}_{y} v|^2
	+ |\nabla\partial^{k-2}_{y} v|^2\right)\, d\Ln
\end{displaymath}
\begin{displaymath}
	\leq \frac{C}{R^4}\int_{B_{4R}^+}
	\left(|\nabla\partial_{y}^{k-1}v|^2 + |\partial_{y}^{k-1}v|^2
	+ |\nabla\partial_{y}^{k-2}v|^2 + |\partial_{y}^{k-2}v|^2\right)\, d\Ln
\end{displaymath}
\begin{displaymath}
	\leq \ldots \leq \frac{C_{\ell,k}}{R^{2\ell}}
	\sum_{j=1}^\ell\int_{B_{2^\ell R}^+}
	\left(|\nabla\partial_{y}^{k-j}v|^2 + |\partial_{y}^{k-j}v|^2\right)\, d\Ln.
\end{displaymath}
By choosing $k=\ell$, (\ref{eq:A_2}) yields
$$\int_{B_R^+}|\partial^k_{y} v|^2\, d\Ln\leq \frac{C_{k}}{R^{2k}}\Ln\left(B_{2^k R}^+\right)\leq C_{k} R^{n-2k}.$$
If $k>2n$, this yields for $R\rightarrow\infty$:
$$\int_{\rel^{n}_+}|\partial^k_{y} v|^2\, d\Ln=0.$$
Hence $\partial^k_{y} v=0$. This implies 
$$\partial_{y_i}\partial^{k-1}_y v=0,\mbox{ for } i=1,\ldots,i-1.$$
Therefore $\partial^{k-1}_{y}v(y,t)$ is independent of $y\in\rel^{n-1}$ and $t\mapsto\partial^{k-1}_{y}v(y,t)$ satisfies the assumptions of Lemma \ref{A_2}. Hence we also have
$$\partial^{k-1}_{y}v=0.$$
 Especially we have
$$\partial_{y_i}\partial^{k-2}_{y}v=0.$$
for all $i=1,\ldots, n-1$.
Hence again $\partial^{k-2}_{y}v(y,t)$ is independent of $y$ and therefore again satisfies the assumptions of Lemma \ref{A_2}.
Iterating this final process yields
$$v=0,$$
which is the desired conclusion.
\proof
{\large \bf Remark:} \\
The proof above works, because both $\Delta^2 $ and $-\Delta$ are monotone operators which are added correctly. If we would destroy this monotonicity by e.g. examining $\Delta^2 + \Delta$, Theorem \ref{half.fourth} is not true anymore. For example a bounded nontrivial solution to
$\Delta^2 v + \Delta v = 0$,
is $(y,t) \mapsto 1 - \cos t$.
\defin

\begin{proposition} \label{half.u}

Let $u:\overline{\rel^{n}_+}\rightarrow \rel$
be a smooth solution of the fourth order boundary-value problem
\begin{equation} \label{half.u.equ}
\begin{array}{c}
	\Delta^2 u - \Delta u = 0
	\quad \mbox{in } \rel^n_+, \\
	u, \partial_n u = 0
	\quad \mbox{on } \rel^{n-1} \times \{0\}. \\
\end{array}
\end{equation}
Furthermore let $u$ for some $\beta \in \rel$ and the corresponding
linear function $u_0: (y,t) \mapsto \beta t$ satisfy
\begin{displaymath}
	u - u_0 \in L^\infty(\rel^n_+).
\end{displaymath}
Then $u$ is one dimensional, that is
\begin{equation} \label{half.u.one}
	u(y,t) = \beta (e^{-t} - 1 + t)
	\quad \mbox{for } y \in \rel^{n-1}, t \geq 0,
\end{equation}
in particular
\begin{equation} \label{half.u.lap}
	\Delta u \equiv \beta
	\quad \mbox{in } \rel^{n-1} \times \{ 0 \}
\end{equation}
and
\begin{equation} \label{half.u.pos}
	u > 0
	\quad \mbox{in } \rel^n_+,
	\quad \mbox{if } \beta > 0.
\end{equation}
\end{proposition}
{\pr Proof:} \\
We put
\begin{displaymath}
	v(y,t) := u(y,t) - \beta (e^{-t} - 1 + t)
	\quad \mbox{for } y \in \rel^{n-1}, t \geq 0,
\end{displaymath}
and see $v \in C^\infty_{loc}(\overline{\rel^n_+})$ and
\begin{displaymath}
\begin{array}{c}
	\Delta^2 v - \Delta v = 0
	\quad \mbox{in } \rel^n_+, \\
	v, \partial_n v = 0
	\quad \mbox{on } \rel^{n-1} \times \{0\}, \\
	v = u - u_0 - \beta (\exp(-.) - 1) \in L^\infty(\rel^n_+). \\
\end{array}
\end{displaymath}
Then the uniqueness in Proposition \ref{half.fourth} gives $v \equiv 0$,
which is (\ref{half.u.one}),
and by direct calculation
\begin{displaymath}
	\Delta u(y,0) = \partial_{tt} (\beta (e^{-t} - 1 + t))|_{t = 0} = \beta
	\quad \mbox{for } y \in \rel^{n-1},
\end{displaymath}
which is (\ref{half.u.lap}).
For $\beta > 0$, we get by the strict convexity of the exponential function
\begin{displaymath}
	u(y,t) = \beta (e^{-t} - 1 + t) > 0
	\quad \mbox{for } y \in \rel^{n-1}, t > 0,
\end{displaymath}
which is (\ref{half.u.pos}).
\proof

%%%%%
%\input{app-w}
%\input{app-h}

\end{appendix}

%%%%%

%%%%%
\phantomsection
\addcontentsline{toc}{section}{References}% Provides a hyperlink in some pdf readers to the bibliography in their content description 

%%%%%


\begin{thebibliography}{99}

\bibitem[ADN59]{adn1}
	{Agmon, S., Douglis, A., Nirenberg, L.,} 
	{Estimates near the boundary for solutions of
	elliptic partial differential equations
	satisfying general boundary conditions I},
	{Communications on Pure and Applied Mathematics},
	{\bf {12}}, {(1959)},
	{pp. 623--727},
	{\doi{10.1002/cpa.3160120405}}.
\bibitem[ADN64]{adn2}
	{Agmon, S., Douglis, A., Nirenberg, L.,} 
	{Estimates near the boundary for solutions of
	elliptic partial differential equations
	satisfying general boundary conditions II},
	{Communications on Pure and Applied Mathematics},
	{\bf {17}}, {(1964)},
	{pp. 35--92},
	{\doi{10.1002/cpa.3160170104}}.
\bibitem[Al1899]{Almansi1899-1}
	{Almansi, E.,} 
	{Sull'integrazione dell'equazine differenziale $\Delta^2\Delta^2=0$},
	{Rom. Acc. L. Rend.},
	{\bf {81}},
	{5}, {(1899)},
	{pp. 104--107}.
\bibitem[AnBuoFr19]{AntunesBuosoFreitag2019}
	{Antunes, P.R., Buoso, D., Freitas, P.,}  
	{On the behavior of clamped plates under large compression},
	{SIAM J. Appl. Math.},
	{\bf {79}}, {(2019)},
	{pp. 1872--1891},
	{\doi{10.1137/19M1249606}}.
\bibitem[Bi1933]{Bickley}
	{Bickley, W.G.,} 
	{Deflexions and vibrations of a circular elastic plate
	under tension},
	{Phil. Mag.},
	{(7)},
	{\bf {15}}, {(1933)},
	{pp. 776--797},
	{\doi{10.1080/14786443309462222}}.
\bibitem[Bo05]{BoggioGreenFunction}
	{Boggio, T.,} 
	{Sulle funzioni di Green d'ordine},
	{m. Rend. Circ. Mat. Palermo},
	{\bf {20}}, {(1905)},
	{pp. 97--135},
	{\doi{10.1007/BF03014033}}.
\bibitem[BucGaz11]{BucurGazzola2011}
	{Bucur, D., Gazzola, F.,} 
	{The First Biharmonic Steklov Eigenvalue:
	Positivity Preserving and Shape Optimization},
	{Milan J. Math.},
	{\bf {79}}, {(2011)},
	{pp. 247--258},
	{\doi{10.1007/s00032-011-0143-x}}.
\bibitem[BuoKe21]{BuosoKennedy2021}
	{Buoso, D., Kennedy, J.,} 
	{The Bilaplacian with Robin boundary conditions}, {(2021)},
	{arXiv:2105.11249}.
\bibitem[CaTa20]{CassaniTarsia}
	{Cassani, D., Tarsia, A.,} 
	{Maximum Principle for Higher Order Operators in General Domains}, {(2020)},
	{arXiv:2011.01091v1}.
\bibitem[DaSw04]{DallAcquaSweers2004}
	{Dall'Acqua, A., Sweers, G.,} 
	{On domains for which the clamped plate system
	is positivity preserving},
	{Differential Equations and Inverse Problems, ed. by Carlos
	Conca, Raul Manasevich, Gunter Uhlmann and Michael Vogelius, AMS}, {(2004)},
	{\doi{10.1090/conm/362/06609}}.
\bibitem[DN55]{DouglisNirenbergInteriorSchauder}
	{Douglis, A., Nirenberg, L.,} 
	{Interior Estimates for Elliptic Systems
	of Partial Differential Equations},
	{Communications on Pure and Applied Mathematics},
	{\bf {VIII}}, {(1955)},
	{pp. 503--538},
	{\doi{10.1002/cpa.3160080406}}.
\bibitem[Du48]{DuffinStrip}
	{Duffin, R.J.,} 
	{On a question of Hadamard concerning super-biharmonic functions},
	{Journal of Mathematics and Physics},
	{\bf {27}}, {(1948)},
	{pp. 253--258},
	{\doi{10.1002/sapm1948271253}}.
\bibitem[Fr79]{LFrank} L.S. Frank, Coercive singular perturbations I: a-priori
estimates, Ann. Mat. Pura Appl. IV 119, (1979), pp. 41--113.
{\doi{10.1007/BF02413170}}
\bibitem[Gar51]{GarabedianEllipse}
	{Garabedian, P.,R.} 
	{A partial differential equation arising in conformal mapping},
	{Pacific J. Math},
	{\bf {1}}, {(1951)},
	{pp. 485--524}.
\bibitem[GazGrSw]{PolyHarmBoundValue}
	{Gazzola, F., Grunau, H.-Ch., Sweers, G.,} 
	{Polyharmonic Boundary Value Problems}, {(2010)},
	{Springer}, 
	{Lecture Notes in Mathematics 1991},
	{\doi{10.1007/978-3-642-12245-3}}.
\bibitem[GT]{gil.tru}
	{Gilbarg, D., Trudinger, N.S.,} 
	{Elliptic Partial Differential Equations of Second Order}, {(1998)},
	{Springer, 3. Edition},
	{Berlin - Heidelberg - New York - Tokyo},
	{\doi{10.1007/978-3-642-61798-0}}.
\bibitem[Gre68]{WGreenlee} W. Greenlee, Rate of convergence in singular perturbations,
Annales de l'Institut Fourier, Tome 18, (1968), pp. 135--191.
{\doi{10.5802/aif.296}}
\bibitem[Gr21]{GrunauGreenEstimate2021}
	{Grunau, H.-Ch.,} 
	{Optimal estimates from below for green functions of higher order
	elliptic operators with variable leading coefficients}, 
	{Arch. Math.}, {(2021)},
	{\doi{10.1007/s00013-021-01597-x}}.
\bibitem[Gr02]{GrunauOneDimensionalCase2002}
	{Grunau, H.-Ch.,} 
	{Positivity, change of sign and buckling eigenvalues in a one-dimensional fourth order model problem},
	{Adv. Differential Equations},
	{\bf {7}}, {(2002)},
	{pp. 177--196}.
\bibitem[GrRo10]{GrunauRobert2010}
	{Grunau, H.-Ch., Robert, F.,} 
	{Positivity and Almost Positivity of Biharmonic Green’s Functions
	under Dirichlet Boundary Conditions},
	{Archive for Rational Mechanics and Analysis},
	{\bf {195}}, {(2010)},
	{pp. 865--898},
	{\doi{10.1007/s00205-009-0230-0}}.
\bibitem[GrRoSw11]{GrunauRobertSweers2011}
	{Grunau, H.-Ch., Robert, F., Sweers, G.,} 
	{Optimal estimates from below for biharmonic Green functions},
	{Proc. Amer. Math. Society},
	{\bf {139}}, {(2011)},
	{pp. 2151--2161},
	{\doi{10.1090/S0002-9939-2010-10740-2}}.
\bibitem[GrRomSw20]{GrunauRomaniSweers2020}
	{Grunau, H.-Ch., Romani, G., Sweers, G.,} 
	{Differences between fundamental solutions of general higher order
	elliptic operators and of products of second order operators}, 
	{Math. Ann.}, {(2020)},
	{\doi{10.1007/s00208-020-02015-3}}.
\bibitem[GrSw14a]{GrunauSweersSignChangeUniform1}
	{Grunau, H.-Ch., Sweers, G.,} 
	{A clamped plate with a uniform weight may change sign},
	{Discrete \& Continuous Dynamical Systems - S},
	{\bf {7}},
	{(4)}, {(2014)},
	{pp. 761--766},
	{\doi{10.3934/dcdss.2014.7.761}}.
\bibitem[GrSw14b]{GrunauSweersSignChangeUniform2}
	{Grunau, H.-Ch., Sweers, G.,} 
	{In any dimension a "clamped plate" with a uniform weight may change sign},
	{Nonlinear Analysis A: T.M.A},
	{\bf {97}}, {(2014)},
	{pp. 119--124},
	{\doi{10.1016/j.na.2013.11.017}}.
\bibitem[Gue81]{MGuegnon} M. Gueugnon, Perturbations singulieres dans les espaces $L^{p}$
pour l'op\'erateur $\varepsilon \Delta ^{2}-\Delta +I$, C.R. Acad. Sci. Paris
\ 293, (1981), pp. 129--131.
\bibitem[H1908a]{Hadamard1}
	{Hadamard, J.,} 
	{M\'{e}moire sur le probl\`{e}me d'analyse relatif \`{a}
	l'\'{e}quibilibre des plaques \'{e}lastiques encastr\'{e}es},
	{{\OE}uvres de Jaques Hadamard, Tome II, CNRS Paris, (1968), pp. 515--641},
	{Reprint of: M\'{e}moire pr\'{e}sent\'{e}s par divers savant a
	l'Acad\'{e}mie des Sciences (2), (1908), 33:1-128}.
\bibitem[H1908b]{Hadamard2}
	{Hadamard, J.,} 
	{Sur certains cas int\'{e}ressants du probl\`{e}me biharmonique},
	{{\OE}uvres de Jaques Hadamard, Tome III, CNRS Paris, (1968), pp. 1297--1299},
	{Reprint of: Atti IV Congr. Intern. Mat. Rome, (1908), pp. 12--14}.
\bibitem[Hi19]{HimmelBachelorThesis}
	{Himmel, B.,} 
	{Die Boggio-Formel f\"ur polyharmonische Dirichlet-Probleme},
	{Bachelor thesis: Advisor H.-Ch. Grunau},
	{Otto-von-Guericke Universit\"at Magdeburg}, {(2019)},
	{\href{http://www-ian.math.uni-magdeburg.de/home/grunau/papers/Himmel\_Bachelor.pdf}{http://www-ian.math.uni-magdeburg.de/home/grunau/papers/Himmel\_Bachelor.pdf}}.
\bibitem[Naj88]{Najman} Najman, B. Singular perturbations of elliptic boundary value
problems in $L^{p}$. Glas. Mat. Ser. III 23(43), (1988), pp. 259--290.
\bibitem[Pu15]{Pulst}
	{Pulst, L.,} 
	{Dominance of Positivity of the Green’s Function associated to
	a Perturbed Polyharmonic Dirichlet Boundary Value Problem
	by Pointwise Estimates},
	{PhD thesis: Advisor H.-Ch. Grunau},
	{Otto-von-Guericke Universit\"at Magdeburg}, {(2015)},
	{\doi{10.25673/4208}}.

\end{thebibliography}
\end{document}